\documentclass[pdftex,12pt]{amsart}
\usepackage[dvips]{graphicx,color}
\usepackage[pdftex,colorlinks=true, pdfstartview=FitV, linkcolor=blue, citecolor=blue, urlcolor=blue]{hyperref}
\usepackage{amsthm}
\usepackage{amsopn}
\usepackage{amssymb}
\usepackage{amsmath}
\usepackage{latexsym}
\usepackage{graphics}
\usepackage[utf8]{inputenc}
\usepackage[T1]{fontenc}
\usepackage{amsfonts}
\usepackage{color}

\newtheorem{thm}{Theorem}[section]
\newtheorem{prop}[thm]{Proposition}
\newtheorem{lem}[thm]{Lemma}
\newtheorem{cor}[thm]{Corollary}
\newtheorem{mainthm}{Theorem}
\theoremstyle{definition}
\newtheorem{defn}[thm]{Definition}

\newcommand{\bth}{\begin{thm}}
\renewcommand{\eth}{\end{thm}}
\newcommand{\bpr}{\begin{prop}}
\newcommand{\epr}{\end{prop}}
\newcommand{\ble}{\begin{lem}}
\newcommand{\ele}{\end{lem}}
\newcommand{\bco}{\begin{cor}}
\newcommand{\eco}{\end{cor}}
\newcommand{\bde}{\begin{defn}}
\newcommand{\ede}{\end{defn}}

\newcommand{\pf}{\noindent{\bf Proof}\hspace{7pt}}

\DeclareMathOperator{\Aut}{Aut}
\DeclareMathOperator{\diag}{diag}
\DeclareMathOperator{\SL}{SL}
\DeclareMathOperator{\SU}{SU}
\DeclareMathOperator{\Sp}{Sp}
\DeclareMathOperator{\Stab}{Stab}

\newcommand{\De}{\Delta}
\newcommand{\cA}{{\mathcal A}}
\newcommand{\cL}{{\mathcal L}}
\newcommand{\cT}{{\mathcal T}}
\renewcommand{\bar}{\overline}
\newcommand{\KK}{{\mathbb K}}
\newcommand{\fk}{\mathbb{F}_q}
\newcommand{\fkbar}{\bar{\mathbb{F}}_q}
\newcommand{\tA}{{\widetilde{A}}}
\newcommand{\vep}{\varepsilon}

\newcommand{\mn}{\ \medskip \newline }

\newcommand{\ead}{\email}
\newenvironment{keyword}{Keywords: \keywords}{}

\newcommand{\sep}{\hspace{1ex}}
\begin{document}
\title[expander graphs from curtis-tits groups]{Expander graphs from Curtis-Tits groups}
\author{Rieuwert J. Blok}
\address[Rieuwert J. Blok]{department of mathematics and statistics\\
bowling green state university\\
bowling green, oh 43403\\
u.s.a.}
\ead{blokr@member.ams.org}

\author{Corneliu G. Hoffman}
\address[Corneliu G. Hoffman]{university of birmingham\\
edgbaston, b15 2tt\\
u.k.}
\ead{C.G.Hoffman@bham.ac.uk}

\author{Alina Vdovina}
\address[Alina Vdovina]{School of Mathematics and Statistics\\
Newcastle University\\
Newcastle-upon-Thyne, NE1 7RU\\
U.K.}
\ead{alina.vdovina@ncl.ac.uk}

\thanks{A. Vdovina was partially supported from EPSRC grant EP/F014945/1}

\maketitle

\begin{abstract}
Using the construction of a non-orientable Curtis-Tits group of type $\tilde{A_n}$, we obtain new explicit families of expander graphs of valency $5$ for unitary groups over finite fields.
\end{abstract}



\section{Introduction}

Expanders are sparse graphs with high connectivity properties.
Explicit constructions of expander graphs have   potential applications
in computer science and
is an area of active research.
One of the most significant recent results on expanders is that Cayley
graphs of finite simple groups are expanders, see \cite{KasLubNik06},\cite{BreGreTao}. More precisely there is a $k$ and $\epsilon>0$  such that every non-abelian finite simple group G has a set of $k$ generators for which the Cayley graph $X(G; S)$ is an $\epsilon$ expander. The size of $k$ is estimated around $1000$. 

The present paper is a byproduct of the investigation
 in \cite{BloHof2009b, BloHof2009c} of Curtis-Tits structures
 and the associated groups. A Curtis-Tits (CT) structure over
 $\fk$ with (simply laced) Dynkin diagram $\Gamma$ over a finite 
set $I$ is an amalgam $\cA=\{G_i,G_{i,j}\mid i,j\in I\}$ whose rank-$1$
 groups $G_i$ are isomorphic to $\SL_2(\fk)$,
 where $G_{i,j}=\langle G_i,G_j\rangle$, and in which $G_i$ and $G_j$ commute 
if $\{i,j\}$ is a non-edge in $\Gamma$ and are embedded naturally in
 $G_{i,j}\cong \SL_3(\fk)$ if $\{i,j\}$ is an edge in $\Gamma$. 
It was shown in \cite{BloHof2009b} that such structures are determined up
 to isomorphisms by group homomorphisms from the fundamental group 
of the graph $\Gamma$ to the group $\Aut(\fk)\times \mathbb{Z}_2\le \Aut(\SL_2(\fk))$. Moreover,
 in the case when the diagram is in fact a cycle, all such structures have
 non-collapsing completions, which are described in \cite{BloHof2009c}. 
It turns out that such groups can be described as fixed subgroups of 
certain automorphisms of Kac-Moody groups. 
This is an important point since they will turn out
 to have Kazhdan's property (T) hence they will give rise to families of expanders. 
Many of these groups will be Kac-Moody groups themselves but some will not.
 In particular again in the case where the diagram is a cycle we obtain a new group 
which turns out to be a lattice in $\SL_{2n}(K)$ for some local
 field $K$ and so it will have property (T).
 Moreover, the group in question will have 
finite unitary groups as quotients, giving a more concrete result 
for unitary groups than~\cite{KasLubNik06}. In particular we have
\begin{mainthm}\label{maintheorem1}
For any $n>1$ and prime power $q$, there exists an $\epsilon>0$  such that for any positive integer $s$, there exists a  symmetric
 set $S_{n,q^s}$ of five generators for  $\SU_{2n}(q^s)$ so that the family of Cayley graphs  $X(\SU_{2n}(q^s), S_{n,q^s})$ 
 forms an $\epsilon$-expanding family.
\end{mainthm}

One can view this result as a generalization of results that derive expander graphs from Kac-Moody groups.
Indeed it is known c.f.~\cite{DymJan2002} that if 
 $q>\frac{(1764)^n}{25}$, then the automorphism group $G(\fk)$ of a Moufang twin-building over $\fk$ has Property (T) and so by Margulis' theorem (Theorem~\ref{margulis}) if $G(\fk)$ admits infinitely many finite quotients, their Cayley graphs form a family of expanders. 
 In fact, by~\cite{CaRe2006}, non-affine Kac-Moody groups of rank $n<q$ over $\fk$ are almost simple, so they don't have finite quotients. Therefore the above result only applies to locally finite Moufang twin-buildings of affine type.
In this case, it is known that if $G$ is a connected almost $\KK$-simple  algebraic group over a local field $\KK$ with rank $\ge 2$, then $G(\KK)$ has property (T) (cf. Theorem 1.6.1 in~\cite{Bek08}). Again this allows one to create a family of expanders for each group $G(\KK)$ and each characteristic.

Our methods have been introduced in \cite{BOU00, GraHorMuh2011, GraMuh08} in 
 a slightly more general setting. 
Theorem~\ref {maintheorem1} is weaker than the types of results in 
 \cite{KasLubNik06} and \cite{BreGreTao} 
in the sense that the rank and the characteristic of the groups need to be fixed.  However, our construction is very explicit and the generating set involved is very small compared to theirs.

\paragraph{Acknowledgement}
The second and third authors would like to thank the Isaac Newton Institute in Cambridge were part of this work was done.

\section{The groups}

Let $V$ be a free $\fk[t,t^{-1}]$-module of rank $2n$ with basis $\{e_i,f_i\mid i\in I\}$. Here $I=\{i=1,\ldots,n\}$ and $\fk[t,t^{-1}]$ denotes the ring of commutative Laurent polynomials in the variable $t$ over a finite field $\fk$.
Recall that a 
$\sigma$-sesquilinear form $\beta$ on $V$ is a map $\beta:V\times V \to \fk$ so that $\beta$ is linear in the first coordinate and $\beta(u, \lambda v+ w)=\sigma(\lambda)\beta(u,v)+\beta(u,w)$ for all $u,v,w\in V$ and $\lambda\in \fk$. Such a form is determined by its values on basis elements. Let $\beta$ be such that, for all $i,j\in I$, 
$\beta(e_{i},e_{j})=\beta(f_{i},f_{j})=0, \beta(e_{i},
f_{j})=t\delta_{ij}$ and $\beta(f_{i},e_{j})=\delta_{ij}$ where $\sigma\in \Aut(\fk[t,t^{-1}])$ fixes each element of $\fk$ and interchanges $t$ and $t^{-1}$.
More precisely
\[G^\tau:=\{ g \in \SL_{2n}(\fk[t,t^{-1}])| \forall x,y \in V,  \beta(gx,{g}y)=\beta(x,y)\}\]
In \cite{BloHof2009c} it was proved that $G^\tau$ is a universal ``non-orientable''
Curtis-Tits group.
It turns out that the group $G^\tau$ has some very interesting
natural quotients and that its action on certain Clifford-like
algebras is related to
phenomena in quantum physics.

The aim of this paper is to prove that the group $G^{\tau}$ has Kazhdan's property (T). This implies that the finite quotients of this group will form families of expanders. 
Before doing this we record the following lemma.
\ble\label{lem:5-generation}
The group $G^\tau$ can be generated with a symmetric set $S$ of size at most $5$.
\ele
\pf
Consider the element $s\in \SL_{2n}(\fk[t,t^{-1}])$ transforming the basis above as follows. For each $i=1, \dots n-1$, $e_i^s=e_{i+1}$ and  $f_i^s=f_{i+1}$,  $e_n^s = f_1$ and  $f_n^s=t^{-1}e_{1}$  
It is not too hard to see that $s$ is in fact an element of $G^\tau$. Moreover consider the subgroup 
$$L_0=\{\diag(A, I_{n-2}, {}^tA^{-1}, I_{n-2})\mid A\in \SL_2(\fk)\}.$$
It follows from~\cite{BloHof2009c}  that $G^\tau$ is generated by $s$ and $L_0$. Now, since $\fk$ is finite, $L_0$ is generated by an involution $x$ and another element $y$. Hence we can take $S=\{x,y,y^{-1}, s, s^{-1}\}$. 
\qed
\mn
Let $\fkbar$ denote the algebraic closure of $\fk$.
For any $a\in \fkbar^*$ consider the specialization map
$\epsilon_a\colon\fk[t,t^{-1}]\rightarrow \fkbar$ given by
$\epsilon_a(f)=f(a)$. The map induces a homomorphism $\epsilon_a\colon\SL_{2n}(\fk[t,t^{-1}]) \to \SL_{2n}(\fk(a)) $. In some instances the map commutes with the
automorphism $\sigma$ so that one can define a
map $\epsilon_a\colon G^\tau \to \SL_{2n}(\fkbar) $.
The most important specialization maps are given by $a=\pm1$ or $a=\zeta$, a $(q^s+1)$-st root of $1$ where $s$ is a positive integer.
In case $a=\pm 1$, the automorphism $\sigma$ is trivial and if $q>2$, the specialization maps yield surjections of $G^\tau$  onto $\Sp_{2n}(\fk)$ and $\Omega^{+}_{2n}(\fk)$.
\ble\label{lem:Gtau onto SU}
Suppose that 
$a\in \bar{\mathbb{F}}_q$ is a primitive $(q^{s}+1)$-st root of $1$, for some positive integer $s$.
Then, $\epsilon_a(G^\tau)=\SU_{2n}(q^s)$.
\ele
\pf
Define  $\tilde V= V\otimes_{\fk[t,t^{-1}]}\fk(a)$ and let $\tilde \beta$ be the respective evaluation of $\beta$. We shall also denote by $\bar \lambda$ the image of $\lambda$ under the Galois automorphism given by $a \mapsto a^{-1}$.
Define the transvection map $T_{v}(\lambda): \tilde V \to \tilde V$ by $T_{v}(\lambda)(x)=x+\lambda\tilde \beta (x,v) v$.
Note that the group $\SU_{2n}(q^{s})$ is generated by the set
$$\cT=\{T_{v}(\lambda) | \lambda\in \fk(a)\mbox{ with }\bar\lambda + a\lambda =0, \mbox{ and } v \in \{\tilde e_1, \ldots, \tilde f_{n} \}\}$$ 
since the elements in 
$\{\langle T_{e_{i}}, T_{f_{i}}\rangle\mid i\in I\}$ generate a weak Phan system (see~\cite{BSh} for details).
Therefore if we can lift each map in $\cT$ to $G^{\tau}$, Lemma~\ref{lem:Gtau onto SU} will be proved. We propose that for each $v \in \{\tilde e_1, \ldots, \tilde f_{n} \}$, the lift of $T_{v}(\lambda)$ would be given by a ``transvection'' map $\Phi_{v}(x)=x+F\beta(x, v)v$ where $F\in \mathbb{F}_{q}[t,t^{-1}]$ satisfies $F(a)=\lambda$.
This map is obviously in $\SL_{2n}(\mathbb{F}_{q}[t,t^{-1}])$ so the only thing one needs to check is the fact that it leaves $\beta$ invariant. 
An immediate computation gives 
\begin{align*}\beta(x,y)-\beta(\Phi_{v}(x), \Phi_{v}(y)) &=\sigma( F)\beta(x,v)\sigma({\beta(y,v)})+ F\beta(x,v)\beta(v,y)\\ & = (\sigma(F) +t F)\beta(x,v)\sigma({\beta(y,v) })\end{align*}
and so the sufficient conditions are  $F(a)=\lambda$ and $\sigma(F)+tF=0$.

To find $F$ we shall need the following. Let $f_a\in \mathbb{F}_q[t]$ be the minimal (monic) polynomial for $a$. Then $\sigma(f_a)=t^{-2s}f_a$. Namely, note that $a$ and $a^{-1}$ are conjugate. Moreover if $b$ is another root of $f_a$ then $b$ is a root of $x^{q^s+1}-1$ so it is a power of $a$ and in particular $b^{-1}$ is also a root of $f_a$ and of course $b\ne b^{-1}$ since otherwise $f_a$ will not be irreducible. In conclusion the roots of $f_a$ come in pairs $b, b^{-1}$. This means that $f_a(0)=1$. Now $\sigma( f_a)= f_a(t^{-1})=t^{-2s}f_a'$ where $f_a'$ is a monic irreducible polynomial that has the same roots as $f_a$ so $\sigma(f_a)=t^{-2s}f_a$.

We now find $F$. Pick a polynomial $P\in \fk[t]$ so that $P(a)=\lambda$. Since $\bar\lambda + a\lambda=0$, $a$ is a root of $\sigma(P)+tP$ and so $\sigma(P)+tP=f_aG$ for some $G\in \mathbb{F}_q[t, t^{-1}]$. Applying $\sigma$ shows that $\sigma(f_a)\sigma(G)=t^{-1}f_aG$ and we get the condition $\sigma(G)=t^{2s-1}G$.
Assume $G=\sum _{i=-r}^{l}a_it^i$, the condition above gives that $r= 2s-1+l$ and $a_{-r+i}=a_{l-i}$ for each $i=1, \dots l+r$.

We now propose to find an element $H \in  \mathbb{F}_q[t, t^{-1}]$ so that $\sigma(f_aH)+ tf_aH= f_aG$. Then, $F= P-f_aH$ will have the property that $F(a)=\lambda$ and $\sigma(F)+tF=0$. The condition on $H$ is that $\sigma (H)t^{-2s}+tH=G$. The conditions on $G$ imply that one choice for $H$ is $$H=t^{-l-2s}+ t^{-l-2s+1}+\dots t^{-s-1}+ (a_{-s+1}-1)t^{-s}+ \dots (a_l-1)t^{l-1}.$$
\qed
\mn
It follows from the next result, that, in order to conclude that $G^\tau$ has property (T),  it suffices to show that $G^\tau$ is a lattice in $\SL_{2n}(k((t)))$.
\bth
\begin{enumerate}
\item {\rm (Theorem 1.4.15 in~\cite{Bek08})}
Let $K$ be a local field. The group $\SL_n (K)$ has Property
(T) for any integer $n \ge 3$.
\item  {\rm (Theorem 1.7.1 in~\cite{Bek08})}
If $G$ is a locally compact group and $H$ is a lattice in $G$ then $H$ has property (T) if and only if $G$ does.
\end{enumerate}
\eth
\noindent To do this we use the methods of \cite{BOU00,GraHorMuh2011,GraMuh08}. The more general argument is briefly described in Remark 7.11 in \cite{GraHorMuh2011}.
For convenience we state Lemma 6.22 and 6.23 from \cite{BloHof2009c}:

\ble\label{lem:equal distances}
Suppose that $c_\vep\in \De$ satisfies $\delta_*(c_\vep,c_\vep^\tau)=w$, let $i\in I$ and suppose that
$\pi$ is the $i$-panel on $c_\vep$.
Then,
\begin{enumerate}
\item\label{lem:equal distances part a} There exists a word $u\in W$ such that
$u(u^{-1})^\tau$ is a reduced expression for $w$.
\item\label{lem:equal distances part b} If $l(s_iw)>l(w)$, then all chambers $d_\vep\in \pi-\{c_\vep\}$ except one satisfy
$\delta_*(d_\vep,d_\vep^\tau)=w$.
The remaining chamber $\check{c}_\vep$ satisfies $\delta_*(\check{c}_\vep,(\check{c}_\vep)^\tau)=s_iws_{\tau(i)}$.
\item\label{lem:equal distances part c} If $l(s_i w)<l(w)$, then all chambers $d_\vep\in \pi-\{c_\vep\}$ satisfy
$\delta_*(d_\vep,d_\vep^\tau)=s_i w s_{\tau(i)}$.
\end{enumerate}
\ele

For each $w\in W$ with $w^{-1}=w^\tau$, let $C_w=\{ c \in \Delta_+ | \delta_*(c,c^\tau)=w\}$.
We now have the following strong version of Theorem 6.16 of ~\cite{BloHof2009c}.

\ble\label{lem:transtive on Cw}
The group $G^\tau$ is transitive on $C_w$ for each $w\in W$ with $w^{-1}=w^\tau$.
\ele
\pf
For $w=1$, this is Theorem 6.16 from~\cite{BloHof2009c}.
We now use induction on the length of $w$.
To prove the induction step let $c,d\in C_w$ and let $i\in I$ be such that $l(s_iw)<l(w)$.
Let $c',d'\in \De_+$ be $i$-adjacent to $c$ and $d$ respectively.
Then, by Lemma~\ref{lem:equal distances} part~\ref{lem:equal distances part c}, $c',d'\in C_{s_iws_{\tau(i)}}$ and $l(s_iws_{\tau(i)})<l(w)$.
By induction there is a $g\in G^\tau$ with $g(c')=d'$. 
By Lemma~\ref{lem:equal distances} part~\ref{lem:equal distances part b} $g(c)=d$.
\qed

\bco\label{cor gtau has T}
For  $n>1$,  $G^\tau$ is a non-uniform lattice in $\SL_{2n}(\fk((t)))$ with  property (T).
\eco

\pf 
We apply Lemma 1.4.2 of \cite{BOU00} to conclude that $G^{\tau}$ is a lattice. Since $\det\colon M_n(\fk((t)))\to \fk((t))$ is a continuous map between locally compact Hausdorff spaces, $\SL_n(\fk((t)))$ is  locally compact.
The group $\SL_n(\fk((t)))$ acts cocompactly on $\Delta$ because it is chamber transitive, and it acts properly discontinuously because the residues of $\De$ are finite.
Moreover, $G^\tau$ is discrete because $G^\tau\cap U_\vep=1$, where $U_\vep$ is the unipotent radical of the Borel group of $\SL_{2n}(\fk[t,t])$ for its action on $\De_\vep$ ($\vep=+,1$). 

By Lemma~\ref{lem:transtive on Cw}, $G^\tau$ acts transitively on the sets $C_w$ and these partition $\Delta_+$. For each $u\in W$ pick an element $c_u$ so that $\delta_*(c_u,c_u^\tau)=u(u^{-1})^\tau$ to parametrise the orbits of $G^\tau$ on $\Delta_+$. It now follows from Lemma 1.4.2 of \cite{BOU00} that $G^{\tau}$ is a lattice if and only if the series
$\sum_{u\in W} |\Stab_{G^\tau}(c_w)|^{-1}$ converges. By Lemma~\ref{lem:equal distances} there are exactly $q^{l(u)}$ elements of $C_{1_w}$ at distance $u$ from $c_u$ and for each of these $q^{l(u)}$ chambers $d$, $c_u$ is the unique chamber in $ C_{u({u^-1})^\tau}$ such that $\delta(c,d)=u$.
Therefore the group $\Stab_{G^\tau}(c_w)$ acts transitively on these $q^{l(u)}$ chambers and $|\Stab_{G^\tau}(c_w)|^{-1}\le (q^{-1})^{l(u)}$.  Thus it suffices to show that the Poincar\'e series of $W$, defined as $W(x)=\sum_{u\in W} x^{l(u)}$ converges for $x=q^{-1}$. It follows from a result of Bott~\cite{Bot1956}  that 
 $W(x)=\frac{1-x^{n+1}}{(1-x)^{n+1}}$, which clearly converges for small $x$.
That $G^\tau$ is non-uniform follows from~\cite[1.5.8]{BasLub2001}
 since $\SL_{2n}(\fk((t)))$ is transitive on $\De_+$ and $G^\tau$ has infinitely many orbits.\qed

\section{The expanders}

\bde Let $X=(V,E)$ be a finite $k$-regular graph with $n$ vertices.
we say that $X$ is an $(n,k,c)$ expander if for any subset $A\subset
V$,
$|\partial A|\ge c(1- \frac{|A|}{N})|A|.$
Here $\partial A =\{v \in V \mid d(v, A)=1\}$.
 \ede

\bth\label{margulis}{\rm (Margulis \cite{Mar73})} Let $\Gamma$ be a finitely generated group that
has property (T). Let $\cL$ be a family of finite index normal
subgroups of $\Gamma$ and let $S=S^{-1}$ be a finite symmetric set
of generators for $\Gamma$. Then the family  $\{X(\Gamma/N, S)\mid N
\in \cL\}$ of Cayley graphs of the finite quotients of $\Gamma$ with
respect to the image of $S$ is a family of $(n,k,c)$ expanders for
$n=|\Gamma/N|, k=|S|$ and some fixed $c>0$.
\eth
\smallskip\noindent
Lemma~\ref{lem:Gtau onto SU} therefore has the following consequence.
\bco \label{cor:expanders}
Let $n>1$. If $S$ be a symmetric generating set for $G^\tau$ then the family of Cayley graphs
$\{X(\SU_{2n}(q^s), S)\mid s\ge 1\}$ is an expanding family. \eco

\bpr\label{prop:almost always non-trivial images}
The
image of any non-trivial $g\in G^\tau$ in $\SU_{2n}(q^s)$ is non-trivial for all but finitely many $s$.
\epr
\pf
Suppose that the image of $g\in G^\tau$ in $\SU_{2n}(q^s)$ is trivial for infinitely many $s$.
Then $\epsilon_{a}(g)=I_{2n}\in \SU_{2n}(\fk(a))$ for infinitely many $a$.
In particular, if $g_{ij}(t)$ is any entry of $g$, then  
 $g_{ij}(a)=\delta_{ij}$ for infinitely many $a$.
But $g_{ij}\in \fk[t,t^{-1}]$, so $g=I_{2n}$. 
\qed
\mn
Finally Lemma~\ref{lem:5-generation}, Lemma~\ref {lem:Gtau onto SU}, and Corollary~\ref{cor:expanders} prove Theorem~\ref{maintheorem1}.

\end{document}